\documentclass[10pt]{article}
\usepackage{mathrsfs}
\usepackage{stmaryrd}
\usepackage{amsfonts,amsmath,amssymb,amscd}
\usepackage[colorlinks,linkcolor=red,anchorcolor=blue,citecolor=green]{hyperref}
\usepackage{shadow}
\usepackage{graphics,graphicx}
\usepackage{color}
\usepackage{enumerate}
\usepackage{longtable}
\usepackage{appendix}
\allowdisplaybreaks

\newtheorem{thm}{Theorem}[section]
\newtheorem{lem}[thm]{Lemma}

\newtheorem{defi}[thm]{Definition}
\newtheorem{remark}[thm]{Remark}
\newtheorem{example}[thm]{Example}
\newtheorem{problem}[thm]{Problem}
\newtheorem{conj}[thm]{Conjecture}

\newtheorem{ass}[thm]{Assumption}
\def\qed{\hfill \rule{4pt}{7pt}}

\def\pf{\noindent {\it{Proof.}\hskip 2pt}}

\begin{document}
\begin{center}
{\large\bf Locality of percolation critical probabilities: uniformly nonamenable case}
\end{center}
\begin{center}
Song He, Xiang Kai-Nan and Zhu Song-Chao-Hao
\vskip 1mm
\footnotesize{School of Mathematical Sciences, LPMC, Nankai University}\\
\footnotesize{Tianjin City, 300071, P. R. China}\\
\footnotesize{Emails: songhe@mail.nankai.edu.cn} (Song)\\
\footnotesize{~~~~~~~~~~~~kainanxiang@nankai.edu.cn} (Xiang)\\
\footnotesize{~~~~~~~~~~~~zsch61@qq.com} (Zhu)
\footnote{The project is supported partially by CNNSF (No. 11271204).

   {\it \ \ MSC2010 subject classifications}. 60K35, 60C05, 82B43.

   {\it \ \ Key words and phrases}. Percolation critical probability, locality, nonamenable transitive graph.}
\end{center}

\begin{abstract}
Let $\{G_n\}_{n=1}^{\infty}$ be a sequence of transitive infinite
connected graphs with $\sup\limits_{n\geq 1} p_c(G_n) < 1,$ where each $p_c(G_n)$ is bond percolation critical probability on $G_n.$
Schramm (2008) conjectured that if $G_n$ converges locally to a transitive infinite
connected graph $G,$ then $p_c(G_n) \rightarrow p_c(G)$ as $n\rightarrow\infty.$
We prove the conjecture when $G$ satisfies two rough uniformities, and $\{G_n\}_{n=1}^{\infty}$ is uniformly nonamenable.
\end{abstract}

\section{Introduction}

Let $G=(V(G), E(G))$ be a locally finite, connected infinite graph. Denote by $p_c(G)$ the critical probability for Bernoulli bond percolation on $G$, that is,
\[p_c(G)= \inf\left\{p \in [0, 1] :\ \mathbb{P}_p (\exists ~\mbox{an infinite component}) > 0\right\}.\]
Here $\mathbb{P}_{p}=\mathbb{P}_{G,p}$ is the relevant product probability for $p$-Bernoulli bond percolation on $G.$ {\it Notice that throughout this paper all graphs are assumed to be connected even if we don't mention this}.

Let $\mathbb{N}$ be the set of all natural numbers, and $\mathbb{Z}$ (resp. $\mathbb{Z}_{+}$) the set of all (resp. nonnegative) integers.
Denote by $B_G(v, R)$ the ball of radius $R\in\mathbb{Z}_{+}$ in $G$ centered at $v$. To continue, let $Aut(G)$ be the set of all automorphisms of $G.$ When $Aut(G)$ acts transitively on $G,$
$G$ is called transitive. Say $G$ is spherically symmetric about its vertex $o$ if for all pairs of vertices $x,y$
at the same distance from $o,$ there is an automorphism of $G$ fixing $o$ and mapping $x$ to $y.$
And for transitive $G,$ $G$ is spherically symmetric means it is spherically symmetric about some (in fact, any) vertex $o.$

Call a sequence $\{G_n\}_{n=1}^{\infty}$ of transitive graphs converges locally to $G$ if for any integer $R\in\mathbb{Z}_{+},$
there exists $N\in\mathbb{N}$ such that $B_{G_n} (v_n, R)$ and $B_G(v, R)$ are isomorphic as rooted graphs
for all $n \geq N$, where $v$ (resp. $v_n$) is an arbitrary vertex of $G$ (resp. $G_n$). In the setting of general
transitive graphs, there is a natural problem about the locality of percolations.

\begin{problem}\label{problem1}
Let $\{G_n\}_{n=1}^{\infty}$ be a sequence of transitive infinite
graphs which converges locally to a transitive infinite graph $G$. Does $p_c(G_n) \rightarrow p_c(G)$?
In another word, is the value of $p_c$ determined by the local geometry or by the global property of graphs?
\end{problem}

The answer to the problem is negative generally. For example, the following holds:
\[\mbox{As}\ n\rightarrow\infty,\ \mathbb{Z}\times \mathbb{Z}/n\mathbb{Z}\ \mbox{converges locally to}\ \mathbb{Z}^2,\]
but
\[1=p_c(\mathbb{Z}\times \mathbb{Z}/n\mathbb{Z})\nrightarrow p_c(\mathbb{Z}^2)<1.\]
O. Schramm considered that Problem \ref{problem1} should
have a positive answer whenever the previous obstruction is avoided.
The following conjecture is stated in \cite{BNP2011} and attributed to Schramm (2008).
\begin{conj}\label{conj1}
Let $\{G_n\}_{n=1}^{\infty}$ be a sequence of transitive infinite
graphs with $\sup\limits_{n\geq 1} p_c(G_n) < 1$ such that $G_n$ converges locally to a transitive infinite graph $G.$ Then $p_c(G_n) \rightarrow p_c(G)$.
\end{conj}

Whether $\sup\limits_{n\geq 1} p_c(G_n) < 1$ is equivalent or not to $p_c(G_n) < 1$ for
all $n$ is unknown. Besides, $0<p_c(H)<1$ for any nonamenable transitive infinite connected graph $H;$ but
no geometric characterization of the
probabilistic condition $p_c(G) < 1$ has been established so far, which constitutes
part of the difficulty of Conjecture \ref{conj1}.
The conjecture suggests that the percolation critical probability is locally determined.
This contrasts with critical exponents which are believed to be universal and depend
only on global properties of the graph. For background and
further conjectures regarding percolation on infinite graphs can refer to \cite{BS1996, BI2013, RLYP2013, PG2014}.

For Conjecture \ref{conj1}, a related but not directly linked fact is the following: For large enough $d,$ every vertex of $\mathbb{Z}^d$ locally feels like in a $2d$-regular tree $T_{2d};$ Hara and Slade \cite{HTSG1990} proved that
$$p_c\left(\mathbb{Z}^d\right)=1/(2d)+o\left(1/d\right)\ \mbox{comes close to}\ p_c\left(T_{2d}\right)=1/(2d-1).$$

Recently, Delfosse and Z\'{e}mor \cite{DZ2014} studied Bernoulli bond percolation on $m$-regular hyperbolic tillings with
$m\geq 5,$ and derived an upper bound $p_h$ of corresponding $p_c.$
They conjectured $p_c=p_h$ (\cite{DZ2014} Conjecture 7.2), which is a local property for critical probabilities of hyperbolic percolation
in a sense close to above Conjecture \ref{conj1}.\\

The locality of the value of $p_c$ is a natural and important question which attracts a lot of attention.  For example $p_c$ for slabs $\mathbb{Z}^d\times \{1, \cdots, k\}^{\ell}$ (with $d \geq 2$ and $\ell\geq 1$) approaches
$p_c(\mathbb{Z}^{d+\ell})$. Grimmett and Marstrand \cite{GM1990} resettled this question
positively by a block construction and renormalization argument. Though the just mentioned Grimmett-Marstrand's result is not in the setting of transitive graphs, as remarked in \cite{BNP2011}, from it, we see easily that as $k\rightarrow\infty,$
$$p_c\left(\mathbb{Z}^d\times (\mathbb{Z}/k\mathbb{Z})^{\ell}\right)\rightarrow p_c\left(\mathbb{Z}^{d+\ell}\right).$$
In \cite{BS1996}, Benjamini and Schramm studied the percolation properties in general setting of transitive graphs.
In \cite{BNP2011}, Benjamini, Nachmias and Peres showed that the percolation is local geometry for non-amenable graphs with large girth;
namely, Conjecture \ref{conj1} holds when $G$ is a $d$-regular tree $T_d$ with $d\geq 3$ and $\{G_n\}_{n=1}^{\infty}$ is uniformly nonamenable.
Here uniform nonamenability of $\{G_n\}_{n=1}^{\infty}$ and spherically symmetrical tree-structure of $T_d$ play an important role.
Martineau and Tassion \cite{MT2013} proved Conjecture \ref{conj1} in case of Cayley graphs of Abel groups with rank no less than 2 by
making the Grimmett-Marstrand argument more robust. Note in \cite{GM1990} and \cite{MT2013}, the graphs are assumed to be amenable. And the conjecture is still open even for uniformly nonamenable sequences $\{G_n\}_{n=1}^{\infty}$ (\cite{BNP2011}).

As mentioned in Pete \cite{PG2014} Section 14.2, Conjecture \ref{conj1} appears to be quite hard. Generally, in \cite{PG2014} Section 14.2,
Pete proved one direction of the Schramm conjecture that
$$p_c(G)\leq\liminf\limits_{n\rightarrow\infty}p_c(G_n)    \eqno{(1.1)}$$
by the following result: For Bernoulli bond percolation on any infinite transitive graph $H$ with that $p_c(H)<1-\delta$ for some $\delta\in (0,1),$
there is a positive constant $c(\delta)$ only depending on $\delta$ such that
$$\theta_H(p)-\theta_H(p_c(H))\geq c(\delta)(p-p_c(H)),\ \forall p\in (p_c(H),1];$$
where $\theta_H(p)=\mathbb{P}_{H,p}[\mbox{cluster of}\ o\ \mbox{is infinite}]$ and $o$ is a fixed vertex of $H.$ For a discuss of the other direction
of the conjecture, refer to \cite{PG2014} Section 14.2.\\

To state our main result, we introduce the following definition of ``{\it quasi-spherically symmetric}".
\begin{defi}
{\bf (i)} Call a graph $H$ is quasi-spherically symmetric about its vertex $o$ if there exists a $k\in\mathbb{N}$ such that
for any $n\in\mathbb{N},$ $\partial B_H(o,n):=B_H(o,n)\setminus B_H(o,n-1)$ can be divided into at most $k$ disjointed subsets
with that for any two vertices $x$ and $y$ in the same subset, there is an automorphism of $H$ fixing $o$ and mapping $x$ to $y.$

{\bf (ii)} And for transitive $H,$ $H$ is quasi-spherically symmetric means it is quasi-spherically symmetric about some (in fact, any) vertex $o.$
\end{defi}

Clearly, for a graph $H,$ spherical symmetry implies quasi-spherical symmetry.
Note any $d$-regular tree $T_d$ with $d\geq 2$ is a spherically symmetric Cayley graph and any Cayley graph is symmetric,
and Cayley graphs are rarely spherically symmetric. However, there still are some spherically symmetric Cayley graphs rather than regular trees; and there are
some quasi-spherically symmetric
transitive graphs which are not spherically symmetric. See Section 2.\\

Given an infinite transitive graph $G=(V(G), E(G))$ and a simple random walk (SRW) $\{X_n\}_{n=0}^\infty$ on it.
For any $n\in\mathbb{N}, u\in V(G), x\in \partial B_{G}(u,n),p\in (p_c(G),1),$
let $\sigma_n=\inf\{t~|X_t\notin B_G(X_0,n)\}$ and
$$
\mu_n(x)=\mathbb{P}[X_{\sigma_n}=x\ \vert\ X_0=u],\
b_n(x,u)=\mathbb{P}_{G,p}\left[x\ \mbox{is connected to}\ u\ \mbox{in}\ B_{G}(u,n)\right]. \eqno{(1.2)}
$$
Notice $\mu_n(\cdot)$ and $b_n(\cdot,u)$ depend on $u,$ while the following assumption is independent of $u.$
\begin{ass}\label{Assumption1}
(Rough uniformity)
{\bf (i)} There exist $0<c, C<\infty$ such that
\[c\leq \frac{\mu_n(x)}{\mu_n(y)}\leq C,\ \forall n\in\mathbb{N},\ \forall x,y \in \partial B_G(u,n).\]

{\bf (ii)} There are positive constants $c_1(p)$ and $c_2(p)$ depending on $p$ satisfying
$$c_1(p)\leq b_n(x,u)/b_n(y,u)\leq c_2(p),\ \forall n\in\mathbb{N},\ \forall x,y \in \partial B_G(u,n).$$
\end{ass}
\vskip 2mm

Note that Conjecture \ref{conj1} is open even if we assume $\{G_n\}_{n=1}^{\infty}$ is uniformly nonamenable
(\cite{BNP2011} p.202 lines -2 and -1) and appears to be quite hard (\cite{PG2014} Section 14.2). In this paper, we prove

\begin{thm}\label{thm1}
Let $\{G_n\}_{n=1}^{\infty}$ be a sequence of uniformly nonamenable transitive infinite
graphs converging locally to a transitive infinite graph $G.$ If $G$ satisfies Assumption \ref{Assumption1}, then as $n\rightarrow\infty,$
\[p_c(G_n) \rightarrow p_c(G).\]
\end{thm}
\vskip 2mm

By (4.1), local limit of any sequence of uniformly nonamenable transitive infinite
connected graphs is nonamenable. Clearly, Assumption \ref{Assumption1} holds for any spherically symmetric transitive infinite connected graph. Naturally, Assumption \ref{Assumption1} should be true for many quasi-spherically symmetric transitive infinite connected graphs. For spherically or quasi-spherically symmetric and nonamenable transitive infinite connected graphs which are not trees and satisfy Assumption \ref{Assumption1}, see Section 2. Assumption \ref{Assumption1} is a technical condition which should be removed, and we can not remove it in this paper. We conjecture that Assumption \ref{Assumption1} holds for quasi-spherically symmetric transitive graphs.

\begin{remark}
Novel aspects of this paper are as follows.

{\bf (i)} Introduce a notion of quasi-spherically symmetric graphs and find some examples of such nonamenable graphs satisfying Assumption \ref{Assumption1} with uniformly nonamenable local approximations.

{\bf (ii)} Our proof of Theorem \ref{thm1} is inspired by \cite{BNP2011}. Comparing with \cite{BNP2011}, in our setting, the loop erasure of $\{X_t\}_{t=0}^{\sigma_n}$ can not be a uniform random non-backtracking (self-avoiding) path, \cite{BNP2011} Corollary 3.1 does not hold, we introduce notion $(\alpha,A,n)$-nice instead of \cite{BNP2011} $(\alpha,A)$-good. And ingeniously, when considering the bond percolation on $G_n$ in Step 1 of proof for
Lemma \ref{lem4}, we introduce a third family  $\{Z_e(\epsilon_1)\}_{e\in E(G_n)}$ of independent Bernoulli random variables with mean $\epsilon_1$, which is independent of $\{X_e(p)\}_{e\in E(G_n)}$ and $\{Y_e(\epsilon)\}_{e\in E(G_n)}$; and add vertices into $\mathcal{V}_t$ not only according to Case 1 like \cite{BNP2011}, but also according to Case 2 which differs from that of \cite{BNP2011}. And in order that conditional expectation of  $\left\vert\mathcal{V}_t\right\vert$ can be as big as possible, we need to estimate the conditional probability of a vertex being in $\mathcal{V}_t$ has a not too small lower bound specified in (3.4); which is done in Step 3 of proving Lemma \ref{lem4}. Note Assumption \ref{Assumption1} is used to deduce (3.2)-(3.3) in Step 2 of proof for Lemma \ref{lem4}, which does not appear in \cite{BNP2011}. The Steps 2 and 3 are essentially new.
Step 4 in proving Lemma \ref{lem4} is the same as that of \cite{BNP2011}.
\end{remark}

\section{Quasi-spherically symmetric nonamenable transitive graphs}

In this section we will present some examples for nonamenable quasi-spherically symmetric graphs (including spherically symmetric ones) satisfying Assumption \ref{Assumption1} with uniformly nonamenable local approximations. As mentioned before, $d$-regular tree $T_d$ is spherically symmetric and nonamenable for any  $d\geq 3$.

\begin{lem}\label{lem0}
Given a nonamenable right Cayley graph $H$ corresponding to a group $\Gamma,$ and a sequence $\{r_n\}_{n=1}^{\infty}\subset\Gamma$ such that
(i) $\vert r_n\vert,$ the distance of $r_n$ from the identity of $\Gamma,$ converges to $\infty$ as $n\rightarrow\infty;$
(ii) additionally when order of $r_n$ is infinite, for some positive constant $c_n,$ $\left\vert r_n^k\right\vert\geq c_n \vert k\vert \vert r_n\vert,\ k\in\mathbb{Z}.$
Let each $H_n$ be the transitive quotient
graph of $H$ given by the following equivalence:
$$x,y\in\Gamma,\ x\sim_n y\Longleftrightarrow x^{-1}y\in\left\{r_n^k:\ k\in \mathbb{Z}\right\}=\langle r_n\rangle;$$
i.e., the $\langle r_n\rangle$-left-coset Cayley graph. Then $\{H_n\}_{n=1}^\infty$ is uniformly nonamenable and locally convergent to $H.$
\end{lem}
\pf For any $h\in H$, let $ord(h)$ be its order. Note
\[p_{H_n}^j\left([x]_n,[x]_n\right)=\displaystyle\sum_{k=-ord(r_n)+1}^{ord(r_n)-1}p_H^j\left(x,xr_n^k\right),\]
where $[x]_n$ is the quotient image of $x$, and $p_{H_n}^j(\cdot,\cdot)$ and $p_H^j(\cdot,\cdot)$ are $j$-step transition probabilities for SRW on $H_n$ and $H$ respectively.
By \cite{PG2014} Chapter $1$ Theorem 1.4, we see that
\[p_H^{j}\left(x,xr_n^k\right)\leq 2\rho^j e^{-\frac{\left|r_n^k\right|^2}{2j}}.\]
Here  $\rho=\rho(H)$ is the spectral radius of $H$.

Therefore, when $ord(r_n)=\infty$, we have that
\begin{align*}
  p_{H_n}^j\left([x]_n,[x]_n\right)&=\sum_{k=-\infty}^{\infty}p_H^j\left(x,xr_n^k\right) \leq p_{H}^j\left(x,x\right)+\sum_{k=1}^{\infty} 4\rho^j e^{-\frac{k^2 c_n^2\left|r_n\right|^2}{2j}}\\
  &\leq 2\rho^j +\sum_{1\leq k\leq j}4\rho^j +\sum_{k\geq j+1}4\rho^j e^{-\frac{kc_n^2\left|r_n\right|^2}{2}}\leq \rho^j\left(2+ 4j+4\sum_{k\geq j+1}e^{-\frac{kc_n^2\left|r_n\right|^2}{2}}\right)\\
  &\leq \rho^j\left(2+ 4j+4\frac{e^{-\frac{(j+1)c_n^2\left|r_n\right|^2}{2}}}{1-e^{-\frac{c_n^2\left|r_n\right|^2}{2}}}\right).
\end{align*}
This implies $\displaystyle\limsup_{j\rightarrow \infty}\left(p_{H_n}^j\left([x]_n,[x]_n\right)\right)^{\frac{1}{j}}\leq \rho.$
Obviously, $\displaystyle\liminf_{j\rightarrow \infty}\left(p_{H_n}^j\left([x]_n,[x]_n\right)\right)^{\frac{1}{j}}\geq\rho.$
Namely, the spectral radius $\rho(H_n)$ of $H_n$ is  $\rho$.

And when $ord(r_n)<\infty$, we see that
\begin{align*}
  p_{H_n}^j\left([x]_n,[x]_n\right)&=\displaystyle\sum_{k=-ord(r_n)+1}^{ord(r_n)-1}p_H^j\left(x,xr_n^k\right)
  \leq \displaystyle\sum_{k=-ord(r_n)+1}^{ord(r_n)-1}2\rho^j\leq 2(2 ord(r_n)-1)\rho^j,
\end{align*}
which implies $\rho(H_n)\leq \rho$, and further $\rho(H_n)= \rho$.

Clearly $\{H_n\}_{n=1}^{\infty}$ converges locally to $H$ due to $\lim\limits_{n\rightarrow\infty}\vert r_n\vert=\infty.$
\qed\\

Let $P$ be the transition matrix of an SRW on a connected graph $H$ and
$I$ the identity matrix, the bottom $\lambda_1(H)$ of the spectrum of $I-P$ is the largest constant $\lambda_1$ such that for all $f\in \ell^2(H,\deg(\cdot))$, we have
\[\langle f, (I-P)f\rangle \geq \lambda_1 \langle f, f\rangle.\]
Here degree function $\deg(\cdot)$ is viewed as a stationary measure for the SRW on $H,$ and $\langle \cdot, \cdot\rangle$
is the inner product on $\ell^2(H,\deg(\cdot)).$ Notice that when $H$ is transitive, $\ell^2(H,\deg(\cdot))$ can be replaced by $\ell ^2(H).$
Recall Dodziuk \cite{DJ1984} proved that an infinite bounded degree connected graph is nonamenable if and only if $\lambda_1(H)>0.$

\begin{example}\label{exm1}
Modified grandparent graphs
\end{example}

Fix $3\leq d \in \mathbb{N}$. Let $\xi$ be a fixed end of $d$-regular tree $T_d$. For any vertex $x$ of $T_d$, there is a unique ray $x_\xi=\langle x_0=x, x_1,x_2,\cdots\rangle$ starting at $x$ such that $x_\xi$ and $y_\xi$  differ by only finitely many vertices for any vertices $x$ and $y$. Call $x_2$ is the $\xi$-grandparent of $x$. Let $H$ be the graph obtained from $T_d$ by adding the edge $xx_2$ between each $x$ and its $\xi$-grandparent. Note that $H$ is a transitive graph but not a Cayley graph of any group. See \cite{RLYP2013} Section 7.1.
{\it Given any vertex $o$ of $T_d.$} Boundaries $\partial B_H(o,n)$ of balls $B_H(o,n)$ are very complicated: When one walks in $H,$ he can at each step decide to take
(i) a ``grandparent step", which goes further, but only in a specific direction, (ii) or a usual step, which can go in any direction one wants.
This makes that $H$ is not quasi-spherically symmetric. We are grateful that G. Kozma pointed this to us.\\

But the following transitive {\it modified grandparent graph} $G$ is quasi-spherically symmetric and satisfies Assumption \ref{Assumption1}:
Let ${\rm dist}_{T_d}(\cdot,\cdot)$ be the graph distance in $T_d.$ And add new edges between any two vertices $x$ and $y$ of $T_d$ with ${\rm dist}_{T_d}(x,y)=2.$
Denote the obtained graph by $G.$

{\it Indeed}, for any $n\in \mathbb{N}$, let
\[ S_1^n=\{x\in T_d\ |\ {\rm dist}_{T_d}(o,x)=2n-1\}, ~ S_2^n=\{x\in T_d\ |\ {\rm dist}_{T_d}(o,x)=2n\}. \]
 Then
\[\partial B_G(o,n)=S_1^n\cup S_2^n.\]
Obviously, for every pair vertices $x$ and $y$ in $S_1^n$ (resp. $S_2^n$), there is an automorphism fixing $o$ and mapping $x$ to $y$ due to the spherical symmetry of $T_d$. Hence $G$ is quasi-spherically symmetric.

Let $P$ and $P_1$ be the transition matrices of SRWs on $G$ and $T_d$ respectively. Let $P_2=(p_2(x,y))_{x,y}$ be transition matrix on $G$ such that
\begin{align*}
p_2(x,y)=
 \left\{
   \begin{array}{ll}
     0, &~{\rm if}\ {\rm dist}_{T_d}(x,y)\not=2,\\
     \frac{1}{d(d-1)}, &~{\rm if}\ {\rm dist}_{T_d}(x,y)=2.
   \end{array}
 \right.
\end{align*}
Then
\[P=\frac{d}{d^2}P_1+ \frac{d(d-1)}{d^2}P_2=\frac{1}{d}P_1+\left(1-\frac{1}{d}\right)P_2.\]
Note $T_d$ is nonamenabe, $\lambda_1(T_d)>0$. Thus for any  $f\in \ell^2(G)=\ell^2(T_d)$,
\begin{align*}
  \langle f,(I-P)f \rangle &= \frac{1}{d} \langle f,(I-P_1)f\rangle+ \left(1-\frac{1}{d}\right)\langle f,(I-P_2)f\rangle  \\
  &\geq \frac{1}{d} \lambda_1(T_d)\langle f,f\rangle.
\end{align*}
Namely, $G$ is nonamenable.

Since $T_d$ is the right Cayley graph of group $\left\langle a_1,\cdots, a_d\left|a_i^2,1\leq i\leq d\right.\right\rangle$. For any $n\in \mathbb{N}$, let $T_d^{(n)}$ be  the $\left\langle (a_1a_2)^{n+1}\right\rangle =\left\{(a_1a_2)^{k(n+1)}:\ k\in\mathbb{Z}\right\}$ left-coset Cayley graph. Then $\left\{T_d^{(n)}\right\}_{n=1}^\infty$ is uniformly nonamenable by Lemma \ref{lem0}. (There are other uniformly nonamenable local approximations to $T_d$. In fact, for
any group $\Gamma_k$ ($k\geq 2$), Olshanskii and Sapir \cite{OS2009} constructed a uniformly nonamenable sequence which converges locally to its Cayley graph. For $k\geq 4,$ Akhmedow \cite{Ak2005} also gave such a sequence.)

Note each $T_d^{(n)}$ is a quotient graph of $T_d$, and the modified grandparent graph $G$ of $T_d$ induces naturally a modified grandparent graph $G_n$ of $T_d^{(n)}$. Clearly,  $\{G_n\}_{n=1}^\infty$ converges locally to $G$. Similarly to prove $G$ is nonamenable, one can check each $G_n$ is also nonamenable and
\[\lambda_1(G_n)\geq \frac{1}{d} \lambda_1\left(T_d^{(n)}\right)\geq \frac{1}{d}\inf_{k\geq 1}\lambda_1\left(T_d^{(k)}\right)>0.\]
Hence $\{G_n\}_{n=1}^\infty$ is uniformly nonamenable.

Note that $\mu_n(\cdot)$ is a constant function on $S_1^n$ or $S_2^n.$ To prove Assumption \ref{Assumption1}(i), it suffices to check
\[\mu_n(x)\geq \mu_n(y) \geq \frac{1}{d^2}\mu_n(x), ~x\in S_1^n, ~y\in S_2^n,~x\sim y,~n\geq 2. \]

To begin, let
\[\Gamma_{x,z}=\left\{\mbox{path}~ \gamma=\gamma_0\gamma_1\cdots\gamma_k:\ k\in \mathbb{N}, \gamma_0=o, \gamma_k=x, \gamma_{k-1}=z,
      \gamma_0\gamma_1\cdots\gamma_{k-1}\subseteq B_G(o,n-1)\right\},\]
\[\Gamma_{x,w}=\left\{\mbox{path}~ \gamma=\gamma_0\gamma_1\cdots\gamma_k:\ k\in \mathbb{N}, \gamma_0=o, \gamma_k=x, \gamma_{k-1}=w, \gamma_0\gamma_1\cdots\gamma_{k-1}
      \subseteq B_G(o,n-1)\right\},\]
where $z$ and $w$ are parent and grandparent of $x$ respectively (assuming $o$ is the ancestor of all other vertices of $T_d$), see Figure $1.$ And let
\[\Gamma_{y,z}=\left\{\mbox{path}~ \gamma=\gamma_0\gamma_1\cdots\gamma_k:\ k\in \mathbb{N}, \gamma_0=o, \gamma_k=y, \gamma_{k-1}=z, \gamma_0\gamma_1\cdots\gamma_{k-1}
   \subseteq B_G(o,n-1)\right\}.\]
For path $\gamma=\gamma_0\gamma_1\cdots\gamma_k,$ write $\vert \gamma\vert =k.$ Then
\[\mu_n(x)=\sum_{\gamma\in \Gamma_{x,z}}\left(\frac{1}{d^2}\right)^{\left|\gamma \right|}+\sum_{\gamma\in \Gamma_{x,w}}\left(\frac{1}{d^2}\right)^{\left|\gamma
      \right|},\]
\[\mu_n(y)=\sum_{\gamma\in \Gamma_{y,z}}\left(\frac{1}{d^2}\right)^{\left|\gamma \right|}.\]

Since for any $\gamma=\gamma_0\gamma_1\cdots\gamma_k\in \Gamma_{x,z}$, $\widehat{\gamma}=\gamma_0\gamma_1\cdots\gamma_{k-1}y\in\Gamma_{y,z}$; and for any
$\beta=\beta_0\beta_1\cdots\beta_{\ell}\in \Gamma_{x,w}$, $\widetilde{\beta}=\beta_0\beta_1\cdots\beta_{\ell-1}zy\in\Gamma_{y,z}$; we have that
   \[\mu_n(y) \geq \frac{1}{d^2}\mu_n(x).\]
On the other hand, for any $\gamma=\gamma_0\gamma_1\cdots\gamma_k\in \Gamma_{y,z}$, $\overline{\gamma}=\gamma_0\gamma_1\cdots\gamma_{k-1}x\in\Gamma_{x,z}$. So
\[\mu_n(x)\geq \mu_n(y).\]

Notice that $b_n(\cdot,o)$ is a constant function on $S_1^n$ or $S_2^n;$ and for any $x\in S_1^n$ (resp. $S_2^n$), there is an its neighbour $y\in S_2^n$
(resp. $S_1^n$). Then $b_n(x_1,o)\geq p b_n(x_2,o)$ and $b_n(x_2,o)\geq p b_n(x_1,o)$ for any $x_1\in S_1^n$ and $x_2\in S_2^n.$ Hence
\[p\leq \frac{b_n(x,o)}{b_n(y,o)}\leq \frac{1}{p},\ \forall x,y\in\partial B_G(o,n),\ \forall n\in\mathbb{N}.\]
Assumption \ref{Assumption1}(ii) holds.
\begin{figure}[ht]
\centering
\includegraphics[]{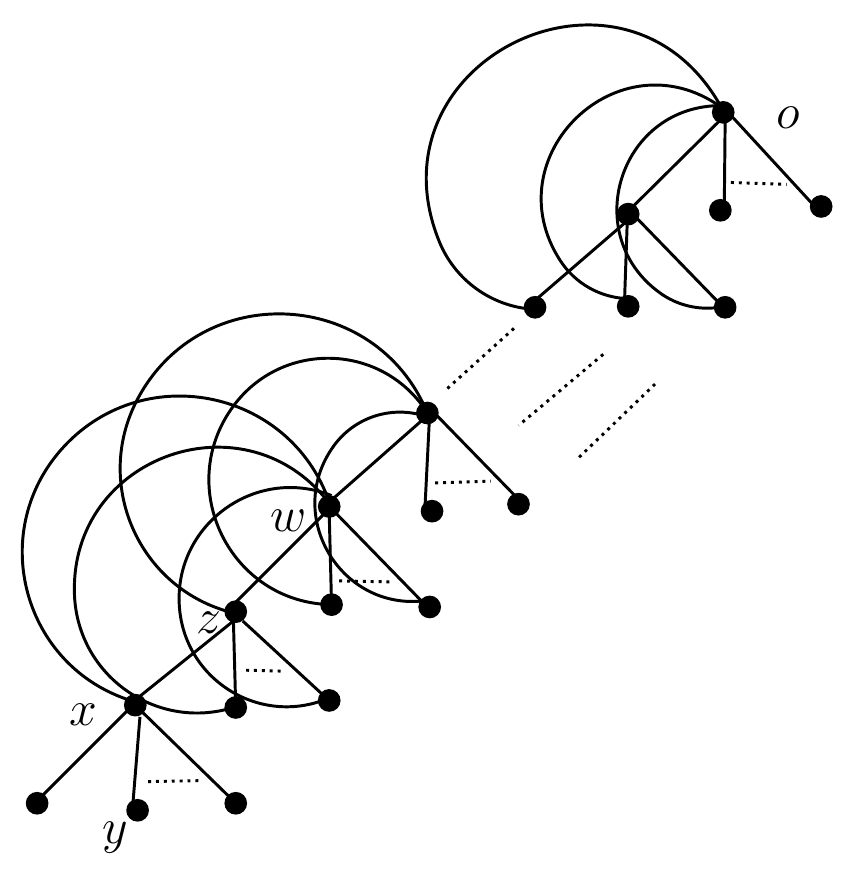}
\caption{Local structure of {\it modified grandparent graph} drawn in plane}
\end{figure}
\qed
\\
\begin{example}\label{exm2}
$\mathbb{Z}_2\ast \mathbb{Z}_3\cong \langle a, b |a^2, b^3\rangle,$ where $\ast$ denotes the free product.
\end{example}

  Cayley graph of above group is shown as  Figure $2$. Given a root $o$, there are only two ways to connect $\partial B_G(o,n)$ with $\partial B_G(o,n-1)$
for any $n\geq 1$. Thus $G$ is quasi-spherically symmetric. Note \cite{WW2000}  Chapter 2 Theorem 10.10:  {\it Let $H=(V(H), E(H))$ be a degree bounded connected infinite graph such that for some $r\in \mathbb{N}$, $H\backslash B_H(x,r)$ has at least three infinite connected components for any $x\in V(H).$ Then $H$ is nonamenable.} By this result, $\mathbb{Z}_2\ast \mathbb{Z}_3$ and $G$  are nonamenable.

Let $G_n$ be the $\left\langle (ab)^n\right\rangle=\left\{(ab)^{kn}:\ k\in\mathbb{Z}\right\}$ left-coset Cayley graph.
Then $\{G_n\}_{n=2}^\infty$ converges locally to $G$; and by  Lemma \ref{lem0}, $\{G_n\}_{n=2}^\infty$ is uniformly nonamenable.

Notice Figure 2. Given any $x, y \in \partial B_G(o,n)$. If $x, y$ are of the same type, then
$$\mu_n(x)=\mu_n(y)\ \mbox{and}\ b_n(x,o)=b_n(y,o).$$
Otherwise, without  loss of generality, assume $x$ is of type $2$ and $y$ is of type $1$. And we can find a type $2$ vertex $z\in \partial B_G(o,n-1)$ and $y\sim z.$ Clearly,
$$\mu_n(y)\leq \mu_{n-1}(z)\ \mbox{and}\ b_n(y,o)\leq b_{n-1}(z,o).$$
Note that
$$\mu_n(x)\geq \frac{1}{3}\mu_{n-1}(z)\ \mbox{and}\ b_n(x,o)\geq pb_{n-1}(z,o).$$
Therefore, we obtain that
$$\mu_n(x)\geq \frac{1}{3}\mu_n(y)\ \mbox{and}\ b_n(x,o)\geq p b_{n}(y,o).$$
Exchanging positions of $x$ and $y$, we have that
$$\mu_n(y)\geq \frac{1}{3}\mu_n(x)\ \mbox{and}\ b_n(y,o)\geq p b_n(x,o).$$
Hence for any $n\in\mathbb{N}$ and $x,y\in\partial B_G(o,n),$
\[\frac{1}{3}\leq\frac{\mu_n(x)}{\mu_n(y)}\leq 3\ \mbox{and}\ p\leq \frac{b_n(x,o)}{b_n(y,o)}\leq \frac{1}{p};\]
Assumption \ref{Assumption1} holds.
\qed\\

\begin{example}\label{exm3}
$H_1\ast H_2 \ast \cdots \ast H_m$ ~($2\leq m\in \mathbb{N}$). Here each $H_i=\langle S_i|R_i\rangle$ ($1\leq i\leq m$) is a nontrivial finite group;
and when $m=2,$ $\max\left\{\vert H_1\vert,\vert H_2\vert\right\}\geq 3.$
\end{example}

Cayley graph $G$ of $H_1\ast H_2 \ast \cdots \ast H_m$ is quasi-spherically symmetric. The reason is as follows: there are finitely many ways for
connecting  $\partial B_G(o,n)$ and $\partial B_G(o,n-1)$ for all $n\geq 1$, where $o$ is a fixed vertex of $G$. From \cite{WW2000}  Chapter 2 Theorem 10.10,  $G$ is nonamenable.

Similarly to Example \ref{exm2}, by Lemma \ref{lem0}, we can construct many uniformly nonamenable graph sequences $\{G_n\}_{n=1}^\infty$  of infinite transitive connected graphs to locally approximate $G$. And also similarly to Example \ref{exm2}, we can check Assumption \ref{Assumption1} holds for $G$.
\begin{figure}[ht]
\centering
\scalebox{0.7}[0.7]{\includegraphics{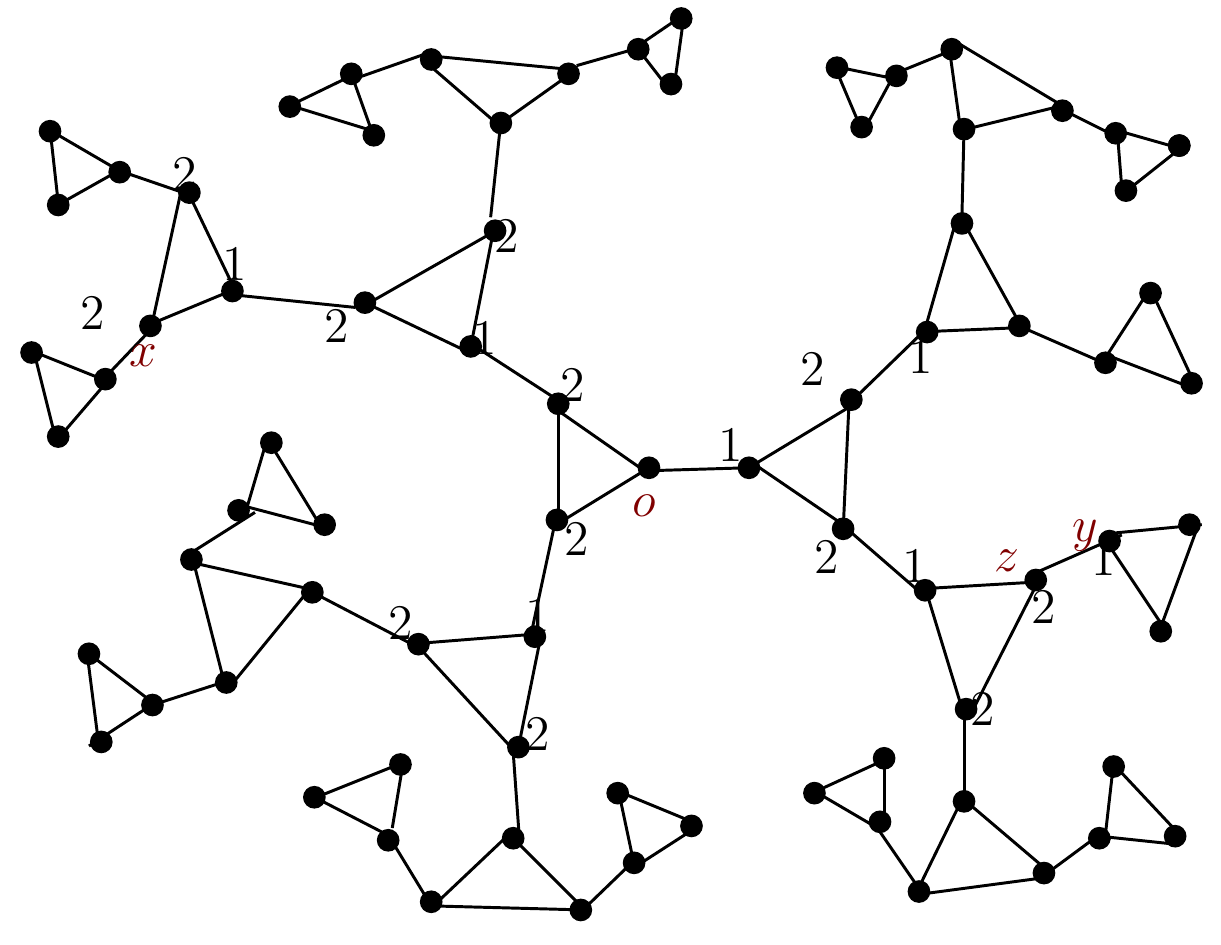}}
\caption{Local structure of Cayley graph for $\mathbb{Z}_2\ast \mathbb{Z}_3$}
\end{figure}
\qed
\begin{example}\label{exm4}
For any composite natural number $d\geq 4,$ there exist $d$-regular spherically symmetric nonamenable infinite Cayley graphs $G$
which are not $d$-regular tree $T_d.$
\end{example}

Since $d$ is a composite number, there exist $2\leq m,k \in \mathbb{N}$ such that $d=mk$. Note complete graph  $\mathbb{K}_{k+1}$ on $k+1$
vertices is a right Cayley graph of some finitely generated group $\langle S_1|R_1\rangle$. Let $G$ be the right Cayley graph of group $H=\langle S_1|R_1\rangle \ast \langle S_1|R_1\rangle \ast \cdots \ast \langle S_1|R_1\rangle$ ($m$ copies). Then $G$ is a spherically symmetric nonamenable infinite Cayley graph. And spherical symmetry is obviously. Indeed, for any $ x\in V(G)$ and $r\geq 1$, $G\backslash B_G(x,r)$ will give birth at least $mk^{r-1}$ infinite connected components; by \cite{WW2000}  Chapter 2 Theorem 10.10, $G$ is nonamenable.

Recall from Lyons and Peres \cite{RLYP2013} Section 3.5, a Cayley graph is spherically symmetric iff it is 2-point homogeneous (i.e., distance transitive \cite{MH1982}) in the sense that there is an automorphism taking $u$ to $w$ and $v$ to $x$ for any vertices $u,v,w,x$ with ${\rm dist}(u,v)={\rm dist}(w,x);$ where ${\rm dist}(\cdot,\cdot)$ is the graph distance. It is still true that an infinite transitive graph $L$ is spherically symmetric iff it is 2-point homogeneous. {\it In fact}, assume transitive $L$ is spherically symmetric and given any  vertices $u,v,w,x$ with ${\rm dist}(u,v)={\rm dist}(w,x).$ Firstly there is an automorphism $\phi_1$ of $L$ taking $u$ to $w$, and then there is an automorphism $\phi_2$ of $L$ fixing $w$ and mapping $\phi_1(v)$ to $x.$ Clearly, automorphism $\phi_2\circ\phi_1$ of $L$ takes $u$ to $w$ and $v$ to $x.$  Recall 2-point homogeneous graphs are characterized by Macpherson \cite{MH1982} Theorem 1.2; from \cite{MH1982} p.63 Definition 1.1 and p.64 Paragraph 1, all these graphs, which are not a tree, are all Cayley graphs just described in this example. We thank R. Lyons for pointing these to us.

Clearly, there are many elements $a$ of $H$ with infinite order satisfying the requirement specified in Lemma \ref{lem0}. Let $G_n$ be the $\left\langle a^n\right\rangle=\left\{a^{kn}:\ k\in\mathbb{Z}\right\}$ left-coset Cayley graph. Then $G_n$ converges locally to $G$ as $n\rightarrow \infty$; and by  Lemma \ref{lem0}, $\{G_n\}_{n=1}^\infty$ is uniformly nonamenable.
\qed\\

\section{Proof of Theorem \ref{thm1}}

Recall the following lemma from \cite{BNP2011}.

\begin{lem}\label{lem2}
Assume $X=\{X_t\}_{t=0}^{\infty}$ is a reversible irreducible Markov chain on a countable state space $V$ with infinite stationary measure $\pi$ and transition matrix $P;$
and the bottom of the spectrum of $I-P$ on $\ell^2(V,\pi)$ is $\lambda_1>0.$ Let $A\subset V$ be nonempty with $\pi(A)<\infty$ and $\pi_A(\cdot)=\pi(A\cap \cdot)/\pi(A).$ Then
\[\mathbb{P}_{\pi_A}(\mbox{$X$ never return to A} )\geq \lambda_1,\]
where $\mathbb{P}_{\pi_A}$ is the law of $X$ with $X_0$ of the law $\pi_A.$
\end{lem}

Let $H=(V(H),E(H))$ be an arbitrary nonamenable infinite transitive connected graph and $X=\{X_t\}_{t=0}^{\infty}$ the SRW on it.  For any $n\in\mathbb{N}$ and $u\in V(H),$ let
\[\tau_n=\inf\{t\geq 1~|X_t\notin B_H(X_1,n)\}.\]
For a subset $A\subset V(H)$ and $\alpha\in (0,1),$ an edge $(x,u)$ is $(\alpha,A,n)$-nice if $x\in A$ and
$$\beta_{(x,u,j)}=\mathbb{P}\left[\left. X_1=u, X_{\tau_j}\notin A\ \right\vert\ X_0=x\right]\geq \alpha,\ \forall 1\leq j\leq n. $$
Call $x\in A$ is $(\alpha,A,n)$-nice if there is a $(\alpha,A,n)$-nice edge $(x,u).$ Let $d$ be the vertex degree of $H.$

\begin{lem}\label{lem3}
Given any finite set $A\subset V(H),$ there are at least $\frac{\lambda_1(H) d}{2}\vert A\vert$ edges $(x,u)$ which are $\left(\frac{\lambda_1(H)}{2},A,n\right)$-nice for any $n\in\mathbb{N}.$
\end{lem}
\pf From Lemma \ref{lem2},
\[\frac{1}{d|A|}\sum_{(x,u):~x\in A}\mathbb{P}\left[\left.X_1=u, X_t\notin A,\forall t\geq 1\ \right\vert\ X_0=x\right]\geq \lambda_1(H).\]
Hence there are at least $\frac{\lambda_1(H)d}{2}\vert A\vert$ edges $(x,u)$ with $x\in A$ such that
$$\mathbb{P}\left[\left.X_1=u, X_t\notin A,\forall t\geq 1\ \right\vert\ X_0=x\right]\geq \lambda_1(H)/2.$$
This implies the lemma.\qed

\begin{lem}\label{lem4}
$\forall \epsilon>0,\epsilon_1>0$, there is an $N\in\mathbb{N}$, for any $n\geq N$, $p_c(G_n)\leq p_c(G)+2\epsilon+\epsilon_1.$
\end{lem}
\pf
Let $\epsilon>0$ and $\epsilon_1>0$ be sufficiently small numbers with $p_c(G)+2\epsilon+\epsilon_1<1.$
Write $p=p_c(G)+\epsilon.$  Let $d$ be the vertex degree of $G.$ Without loss of generality, assume each $G_n$ is of vertex degree $d.$
Put
$$r_n=\sup\left\{r\in \mathbb{Z}_{+}:\ B_{G_n}(v_n,r)\cong B_{G}(v,r)\right\}, \eqno{(3.1)}$$
where $v_n$ (resp. $v$) is an arbitrary vertex of $G_n$ (resp. $G$). Note $r_n$ does not depend on $v_n$ and $v.$
Since $G_n$ converges locally to $G,$ we have $\lim\limits_{n\rightarrow\infty}r_n=\infty.$
For convenience, we also assume vertex sets of each $G_n$ and $G$ are identical.\\

{\bf Step 1.} Consider the bond percolation on $G_n$. For each edge $e$ we consider three independent Bernoulli random variables, $X_e(p),$ $Y_e(\epsilon)$
and $Z_e(\epsilon_1)$ with means $p,$ $\epsilon$ and $\epsilon_1$ respectively. And the family $\{(X_e(p),Y_e(\epsilon),Z_e(\epsilon_1))\}_{e\in E(G_n)}$
is independent. An edge $e$ is open if one of the three variables
$X_e(p),$ $Y_e(\epsilon)$ and $Z_e(\epsilon_1)$ takes the value $1$ and closed otherwise. An edge $e$ is $p$-open if  $X_e(p)=1,$ and $\epsilon$-open if  $Y_e(\epsilon)=1,$ and $\epsilon_1$-open if $Z_e(\epsilon_1)=1.$
So the probability of an edge being closed is $(1-p)(1-\epsilon)(1-\epsilon_1).$ Therefore, for a vertex $v\in V(G_n)$, the open cluster $C(v)$ of $v,$ is dominated by $p_c(G)+2\epsilon+\epsilon_1$-bond percolation. In the following we will prove that with positive probability $|C(v)|=\infty$
for large enough $n.$

We will construct the following process, which produces an increasing sequence $\{A_t\}_{t}$ of connected vertex sets satisfying $A(t)\subseteq C(v)$ for all $t.$ Let $A_0$ be the $p$-cluster of $v.$ Call an edge $\epsilon$-unchecked (resp. $\epsilon_1$-unchecked) if we don't know whether it is $\epsilon$-open or $\epsilon$-closed (resp. $\epsilon_1$-open or $\epsilon_1$-closed). Suppose $A_0$ is finite (otherwise the proof is completed) and all the edges touching $A_0$ are $\epsilon$-unchecked and $\epsilon_1$-unchecked. For $t\geq 1$, let $\mathcal{E}_{t-1}$ be the set of $\epsilon$-unchecked edges $(x,u)$ such that $(x,u)$ is $\left(\frac{\lambda_1(G_n)}{2}, A_{t-1},r_n\right)$-nice. If $\mathcal{E}_{t-1}=\emptyset$, then stop the process. Otherwise, we choose some edge $(x,u)\in \mathcal{E}_{t-1} $  by some order and check whether the edge is $\epsilon$-open or $\epsilon$-closed. If it is closed, then let $A_t=A_{t-1}.$ If not, we consider the $p$-percolation on $B_{G_n}(u,r_n)$ and the $\epsilon_1$-percolation on $B_{G_n}(u,r_n)\cap \partial_E A_{t-1}.$ Here  $\partial_E A_{t-1}$ is the set of  edges with only one end in $A_{t-1}$. Let $\mathcal{V}_t$ be
the set of following vertices in $B_{G_n}(u,r_n)\setminus A_{t-1}:$

\ \ \ \ \ \ \ \ {\bf Case 1.} All vertices of any $p$-percolation path starting from $u$ in $B_{G_n}(u,r_n)$ and

\ \ \ \ \ \ \ \ \ \ \ \ not intersecting with $A_{t-1}.$ Trivially, $u$ is in $\mathcal{V}_t.$

\ \ \ \ \ \ \ \ {\bf Case 2.}  All vertices $y$ of $B_{G_n}(u,r_n)$ such that there is a vertex $z\in B_{G_n}(u,r_n)\setminus A_{t-1}$

\ \ \ \ \ \ \ \ \ \ \ \ adjacent to $A_{t-1}$ by an $\epsilon_1$-open edge $e_z$ satisfying $y$ is connected to $z$ by

\ \ \ \ \ \ \ \ \ \ \ \ a $p$-percolation path in $B_{G_n}(u,r_n)$ avoiding $A_{t-1}.$ Here we don't need to

\ \ \ \ \ \ \ \ \ \ \ \ know whether other edges sharing a vertex with
                        $e_z$ are $\epsilon_1$-open or not. \\
Let
\[A_t=A_{t-1}\cup \mathcal{V}_t.\]\\

{\bf Step 2.} Given edge $(x,u)\in {\cal E}_{t-1}.$
For any $y\in \partial B_{G_n}(u,j)$ with $j\leq r_n,$ let
$$
b_j(y,u,n)=\mathbb{P}_{G_n,p}\left[y\ \mbox{is connected to}\ u\ \mbox{in}\ B_{G_n}(u,j)\right].
$$
Since $\delta:=\mathbb{P}_{G,p}\left[u\ \mbox{is connected to infinity}\right]>0$ and $j\leq r_n,$
we have that
\begin{eqnarray*}
&&\mathbb{P}_{G_n,p}\left[u\ \mbox{is connected to}\ \partial B_{G_n}(u, j)\ \mbox{in}\ B_{G_n}(u,j)\right]\\
&&\ \ =\mathbb{P}_{G,p}\left[u\ \mbox{is connected to}\ \partial B_{G}(u, j)\ \mbox{in}\ B_G(u,j)\right]\geq \delta,
\end{eqnarray*}
and further
$$\sum\limits_{y\in\partial B_{G_n}(u,j)}b_j(y,u,n)\geq \delta.$$
By Assumption \ref{Assumption1}(ii), there is a positive constant $c_3(p)$ depending on $p$ such that for any
$\eta\in (0,1),$ any $j\leq r_n$ and any $A\subseteq \partial B_{G_n}(u,j)$ with $\vert A\vert/\vert \partial B_{G_n}(u,j)\vert\geq \eta,$
$$
\sum\limits_{y\in A}b_j(y,u,n)\geq c_3(p)\eta\delta.\eqno{(3.2)}
$$

By Assumption \ref{Assumption1}(i) and Lemma \ref{lem3},
conditioned on $A_{t-1}$ and $(x,u)\in {\cal E}_{t-1},$ there is a positive constant $c_4$ independent of $A_{t-1}$ and $(x,u)$ such that
for any $j\leq r_n,$
$$\vert \left\{y\in \partial B_{G_n}(u,j):\ y\notin A_{t-1}\right\}\vert \geq c_4\lambda_1(G_n)\vert \partial B_{G_n}(u,j)\vert. \eqno{(3.3)}$$ \\

{\bf Step 3.} Define
\[Z_t=|\{e: \mbox{$e$ is an $\epsilon$-closed and $\epsilon$-checked edge touching $A_t$}\}|.\]
Let
\[\tau =\min \left\{t:\ \vert A_t\vert\leq \frac{2t}{\lambda_1(G_n) d}\right\}.\]
Note that we only check the $\epsilon$-status of one edge at each step, thus $Z_t\leq t$. By Lemma \ref{lem3}, if $|A_t|> \frac{2t}{\lambda_1(G_n) d}$, then there exists at least one $\epsilon$-unchecked edge. Let $\mathcal{F}_t$ be the $\sigma$-algebra generated by the $p,\ \epsilon$ and $\epsilon_1$ statuses of the edges up to time $t$ and $\xi_t=|A_{t+1}|-|A_{t}|.$\\

By the definition of $\mathcal{V}_t,$ we have that for
any vertex $y\in \partial B_{G_n}(u,j)\setminus A_{t-1}$ with $j\leq r_n,$
$$\mathbb{P}[y\in\mathcal{V}_t\ \vert\ {\cal F}_{t-1},\tau >t]\geq \epsilon_1b_j(y,u,n).\eqno{(3.4)}$$
{\it Indeed}, fix $A_{t-1}$ and $(x,u),$ and consider the $p$-bond percolation $\omega(u,n)$ on edges of $B_{G_n}(u,r_n)$ which do not touch
$A_{t-1}.$ Let
$$\partial ^{+}A_{t-1}=\left\{y\in G_n\setminus A_{t-1}:\ \exists z\in A_{t-1},\ y\sim z\right\},\ B_t=\partial^{+}A_{t-1}\cap
   B_{G_n}(u,r_n).$$
 For any $y\in B_{G_n}(u,r_n)\setminus \left(A_{t-1}\cup\{u\}\right),$ let
$D_t(y)=\{y\ \mbox{is connected to}\ B_t\ \mbox{in}\ \omega(u,n)\}$ and
\begin{eqnarray*}
&&D_t^1(y)=\{y\ \mbox{is connected to}\ u\ \mbox{with an open path
  avoiding}\ A_{t-1}\cap B_{G_n}(u,r_n)\ \mbox{in}\ \omega (u,n)\},\\
&&D_t^2(y)=D_t(y)\setminus D_t^1(y).
\end{eqnarray*}
When $D_t^1(y)$ holds, $y$ is in $\mathcal{V}_t.$ And when $D_t^2(y)$ holds, there must be an open path
connecting $y$ with a vertex $z\in B_t\setminus\{u\}$ and avoiding $A_{t-1}$ in $\omega(u,n);$ assume $z\sim w\in A_{t-1}$
and let edge $zw$ be $\epsilon_1$-open, then $y\in\mathcal{V}_t.$ Therefore,
\begin{eqnarray*}
\mathbb{P}[y\in\mathcal{V}_t\ \vert\ {\cal F}_{t-1},\tau >t] &\geq &
   \mathbb{P}\left[\left. D_t^1(y)\ \right\vert\ {\cal F}_{t-1},\tau>t\right]+\epsilon_1\mathbb{P}\left[\left. D_t^2(y)\ \right\vert\
                    {\cal F}_{t-1},\tau>t\right]\\
   &\geq &\epsilon_1\mathbb{P}\left[\left. D_t(y)\ \right\vert\ {\cal F}_{t-1},\tau>t\right].
\end{eqnarray*}

Extend $\omega(u,n)$ to a $p$-bond percolation $\widetilde{\omega}(u,n)$ on $B_{G_n}(u,r_n)$ by letting
edges touching $A_{t-1}\cap B_{G_n}(u,r_n)$ be $p$-open independently and independent of all $X_{\cdot}(p),Y_\cdot(\epsilon)$ and
$Z_\cdot(\epsilon_1).$ Let
\begin{eqnarray*}
&&F_t(y)=\left\{y\ \mbox{is connected to}\ u\ \mbox{in}\ \widetilde{\omega}(u,n)\right\},\\
&&F_t^1(y)=\left\{y\ \mbox{is connected to}\ u\ \mbox{by an open path avoiding}\ A_{t-1}\cap B_{G_n}(u,r_n)\ \mbox{in}\
      \widetilde{\omega}(u,n)\right\},\\
&&F_t^2(y)=F_t(y)\setminus F_t^1(y).
\end{eqnarray*}
Then $F_t^1(y)=D_t^1(y).$

Now assume $F_t^2(y)$ holds. Then there must be an open path $\gamma=(y_0y_1\cdots y_i)$ in $\widetilde{\omega}(u,n)$
such that $y_0=y,\ y_i=u$ and some $y_j\in A_{t-1}\cap B_{G_n}(u,r_n).$ Let
$$j_{*}=\min\left\{1\leq j\leq i:\ y_j\in A_{t-1}\cap B_{G_n}(u,r_n)\right\}.$$
Clearly $1\leq j_{*}<i.$ When $j_{*}>1,$ $(y_0y_1\cdots y_{j_{*}-1})$ is an open path in $\widetilde{\omega}(u,n)$ avoiding $A_{t-1}\cap B_{G_n}(u,r_n);$
hence it is also an open path in $\omega(u,n)$ avoiding $A_{t-1}\cap B_{G_n}(u,r_n).$
Combining with $y_{j_{*}-1}\in B_t$ when $j_{*}>1,$  we see $D_t^2(y)$ holds.
In addition, when $j_{*}=1,$ clearly $D_t^2(y)$ holds. Therefore,
$$F_t^2(y)\subseteq D_t^2(y).$$
And further
$$F_t(y)\subseteq D_t(y).$$
So we have that
\begin{eqnarray*}
\mathbb{P}[y\in\mathcal{V}_t\ \vert\ {\cal F}_{t-1},\tau >t]
  &\geq &\epsilon_1\mathbb{P}\left[\left. D_t(y)\ \right\vert\ {\cal F}_{t-1},\tau>t\right]\\
  &\geq &\epsilon_1\mathbb{P}\left[\left. F_t(y)\ \right\vert\ {\cal F}_{t-1},\tau>t\right]\\
  &=&\epsilon_1\mathbb{P}^{\widetilde{\omega}(u,n)}[F_t(y)],
\end{eqnarray*}
where $\mathbb{P}^{\widetilde{\omega}(u,n)}$ is the law of $\widetilde{\omega}(u,n),$ and we have used that
$\widetilde{\omega}(u,n)$ is independent of ${\cal F}_{t-1}$ given $\tau>t$ and $(x,u)\in \mathcal{E}_{t-1}.$
This implies (3.4).\\

Now by (3.2)-(3.4),
$$
\mathbb{E}[\xi_{t-1}\ |\ \mathcal{F}_{t-1}, \tau>t] \geq
                  \sum\limits_{j=1}^{r_n}\epsilon_1\sum\limits_{y\in\partial B_{G_n}(u,j)\setminus A_{t-1}}b_j(y,u,n)
  \geq \epsilon_1 c_4\lambda_1(G_n)c_3(p)c_4\lambda_1(G_n)\delta r_n. \eqno{(3.5)}
$$
\vskip 3mm

{\bf Step 4.} Note by the uniform nonamenability of $\{G_k\}_{k=1}^{\infty},$ we have $\inf\limits_{k\geq 1}\lambda_1(G_k)>0.$
Hence, for large enough $n,$ by (3.5),
$$\mathbb{E}[\xi_{t-1}\ \vert\ \mathcal{F}_{t-1},\tau>t]\geq 4d^{-1}\lambda_1(G_n)^{-1}.$$

Let $X_i=\sum\limits_{j=0}^i(\xi_j-\mathbb{E}(\xi_j)),\ i\in\mathbb{Z}_{+}.$
Note the randomness of $\xi_{i+1}$ is independent of $\mathcal{F}_{i}.$ Then
\begin{align*}
\mathbb{E}(X_{i+1} |\mathcal{F}_i)=X_i+\mathbb{E}\left(\left.\xi_{i+1}-\mathbb{E}\left(\xi_{i+1}\right) \right| \mathcal{F}_i\right)
   =X_i;
\end{align*}
and $\{X_i\}_i$ is a martingale. Clearly, for any $i,$ $\vert X_{i+1}-X_{i}\vert \leq d(d-1)^{r_n}.$
When $n$ is large enough, by the Azuma-Hoeffding inequality (\cite{AS2004} Chapter $7$),  for any $t>1,$
\begin{align*}
 \mathbb{P}(\tau =t+1 |\vert A_0\vert <\infty)&\leq \mathbb{P}\left(\left.\left\vert A_0\right\vert+\sum_{i=0}^{t}\xi_i \leq \frac{2(t+1)}{d\lambda_1(G_n)}\right\vert \vert A_0\vert <\infty\right)\\
  &\leq \mathbb{P}\left(\left.\sum_{i=0}^{t}\xi_i- \sum_{i=0}^{t}\mathbb{E}(\xi_i)\leq \frac{2(t+1)}{d\lambda_1(G_n)}\right\vert
      \vert A_0\vert<\infty \right)\leq e^{-c_5(t+1)},
\end{align*}
where $c_5:=2\lambda_1(G_n)^{-2}d^{-4}(d-1)^{-2r_n}>0.$ For any $K\in\mathbb{N},$ there is a positive probability with $|A_0|\geq K,$
thus for $K$ large enough, we have
\begin{align*}
  \mathbb{P}(\tau= \infty)&\geq \mathbb{P}(|A_0|\geq K) \mathbb{P}(\tau= \infty ||A_0|\geq K)\\
               &\geq \mathbb{P}(|A_0|\geq K)\left[ 1-\sum_{t\geq \frac{\lambda_1(G_n) d K}{2}}e^{-c_5(t+1)}\right]>0.
\end{align*}
Here we have used the fact that $\tau>\frac{\lambda_1(G_n)d\vert A_0\vert}{2}$ due to $\vert A_t\vert$ is increasing.
Clearly, $\{\tau=\infty\}$ implies $\{|C(v)|=\infty\}.$ Hence, there is a positive probability of an infinite cluster in $p_c(G)+2\epsilon+\epsilon_1$-bond percolation on $G_n.$ That is
$p_c(G_n)\leq p_c(G)+2\epsilon+\epsilon_1.$\qed\\

Similarly to Lemma \ref{lem4}, one can prove

\begin{lem}\label{lem5}
For any $\epsilon>0$ and $\epsilon_1>0,$ when $n$ is large enough, $p_c(G)\leq p_c(G_n)+2\epsilon+\epsilon_1.$
\end{lem}

{\bf By Lemmas \ref{lem4}-\ref{lem5}, or Lemma \ref{lem4} and (1.1), we obtain Theorem \ref{thm1} immediately.\qed}

\begin{remark}
Notice (3.1) and (3.5). From the proof of Theorem \ref{thm1}, the following holds: For any sequence $\{G_n\}_{n=1}^{\infty}$ of nonamenable transitive
infinite connected graphs converging locally to a nonamenable transitive infinite connected graph $G,$ if
$$\mbox{$G$ satisfies Assumption \ref{Assumption1} and $\lim\limits_{n\rightarrow\infty}\lambda_1(G_n)^2r_n=\infty,$}$$
then $\lim\limits_{n\rightarrow\infty}p_c(G_n)=p_c(G).$
\end{remark}

\section{Problems}
Let $\{H_n\}_{n=1}^{\infty}$ be a sequence of transitive infinite
connected graphs converging locally to a transitive infinite connected graph $H.$
Fix vertex $v$ (resp. $v_n$) of $H$ (resp. $H_n$). Consider a copy $\widetilde{H}_n$ of $H_n$ on vertex set $V(H)$
such that
$$B_{\widetilde{H}_n}(v,r_n)=B_{H}(v,r_n)\ \mbox{with}\ r_n=\sup\left\{r\in \mathbb{Z}_{+}:\ B_{H_n}(v_n,r)\cong B_{H}(v,r)\right\}.$$
Since $H_n$ converges locally to $H,$ we have that $\lim\limits_{n\rightarrow\infty}r_n=\infty.$
Let $P$ (resp. $P_n$) be the transition matrix of the SRW on $H$ $\left(\mbox{resp.}\ \widetilde{H}_n\right).$ Then
$P_n$ converges to $P$ pointwisely as $n\rightarrow\infty.$ Notice
$$
\lambda_1(H)=\inf\limits_{f\in\ell ^2(V(H))\setminus\{0\}}\langle f, (I-P)f\rangle/\langle f, f\rangle,\
\lambda_1\left(H_n\right)=\inf\limits_{f\in\ell ^2(V(H))\setminus\{0\}}\langle f, (I-P_n)f\rangle/\langle f, f\rangle.
$$
It is easy to prove
\[\limsup\limits_{n\rightarrow \infty}\lambda_1(H_n)\leq \lambda_1(H). \eqno{(4.1)}\]
Generally $\lim\limits_{n\rightarrow\infty}\lambda_1(H_n)=\lambda_1(H)$ may not hold. A natural question is
\begin{problem}\label{pr1}
Given any transitive infinite connected graph $H.$ {\bf (i)} Is there a nontrivial sequence $\{H_n\}_{n=1}^{\infty}$ of uniformly nonamenable transitive infinite connected graphs
converging locally to $H$ when it is nonamenable? {\bf (ii)} Is there a nontrivial sequence $\{H_n\}_{n=1}^{\infty}$ of amenable transitive infinite connected graphs converging locally to amenable $H$ with $p_c(H)<1$ satisfying that $\sup\limits_{n\geq 1}p_c(H_n)<1?$
\end{problem}

See Section 2 for partial positive answers to Problem \ref{pr1}(i). Recall that
every finitely generated infinite group with free group $F_2$ as a subgroup is nonamenable; and the Burnside group
$$B(m,n)=\left\langle g_1,\cdots, g_m\ \left\vert\ g_i^n=1,\ \forall i\right.\right\rangle$$
is nonamenable and does not contain $F_2$ as a subgroup for any $m\geq 2$ and odd $n\geq 665;$
and grandparent graphs are not isomorphic to any Cayley graph, and are transitive and nonamenable.
And as a transitive graph which is not quasi-isometric to any Cayley graph (\cite{EFW2013}), the Diestel-Leader graph $DL(k,\ell)$ is nonamenable iff $k\neq \ell.$
Therefore, there are plenty of typical examples for studying Problem \ref{pr1}(i).

Notice any quotient graph of an amenable graph is amenable.
And for a finitely generated infinite group, there may not be an element of infinite-order generally.
When $H$ is a Cayley graph of an amenable finitely generated infinite group $\Gamma=\langle S\vert R\rangle$
with an element $r$ of infinite order,
then Cayley graph $H_n$ of $\left\langle S\left\vert R, r^n\right.\right\rangle$ converges locally to $H;$
and to answer affirmatively Problem \ref{pr1}(ii), we need $(a)$ each $H_n$ is infinite, and $(b)$ $\sup\limits_{n\geq 1}p_c(H_n)<1.$
Clearly $(a)$ is easy usually when $p_c(H)<1.$ But $(b)$ might be difficult in some cases.
In addition, we point out that solvable groups and groups with subexponential
growth are amenable.

For a sequence $\{H_n\}_{n=1}^{\infty}$ of uniformly nonamenable transitive infinite connected graphs,
clearly $\sup\limits_{n\geq 1}p_c(H_n)<1.$ Conversely, for a sequence $\{H_n\}_{n=1}^{\infty}$ of nonamenable transitive infinite connected graphs
converging locally to a nonamenable transitive infinite connected graph with $\sup\limits_{n\geq 1}p_c(H_n)<1,$
is the sequence uniformly nonamenable? Maybe the answer is negative.\\

In addition, the following problem is fundamental:

\begin{problem}\label{pr2}
What kind of transitive infinite connected graphs satisfy Assumption \ref{Assumption1}(i) or Assumption \ref{Assumption1}(ii)
or Assumption \ref{Assumption1}?
\end{problem}

Assumption \ref{Assumption1}(i) means harmonic measures (exiting distributions) $\mu_n(\cdot)$ on $\partial B_G(u,n)$
are roughly uniform. Recall Pete \cite{PG2014} Section 9.5 proposed the following rough uniformity problem for finitely generated
groups: {\it Does every finitely generated group have a generating set in which harmonic measures $\mu_n(\cdot)$ on $\partial B_G(u,n)$
are roughly uniform in the sense that there exist constants $0<c,C<\infty$ such that for any $n,$ there is $U_n\subseteq\partial B_G(u,n)$
satisfying
$$\mu_n(U_n)>c,\ c<\frac{\mu_n(x)}{\mu_n(y)}<C\ \mbox{for all}\ x,y\in U_n?$$}

Note Graphs specified in examples of Section 2 are hyperbolic and have an infinite hyperbolic boundary. Recall every
hyperbolic transitive graph with infinite hyperbolic boundary is nonamenable (\cite{WW2000} Chapter 4 Section 22). We conjecture
that Assumption \ref{Assumption1} holds for quasi-spherically symmetric transitive graphs. As for hyperbolic transitive graphs of infinite hyperbolic boundary
(e.g., $d$-regular infinite hyperbolic tiling with $d\geq 5$), R. Lyons thinks this is not true for either Assumption \ref{Assumption1}(i) or Assumption \ref{Assumption1}(ii), and the problem is that differences build up multiplicatively
from $u$ to the boundary of the ball. For the same reason, he doubt Pete's question has a positive answer, and like wise for the following
Problem  \ref{pr4}. However for $d$-regular infinite hyperbolic tiling with $d\geq 5,$ we believe that we can prove Problem \ref{pr4}(ii) has a positive answer, which we plan to do in another paper.

I. Benjamini suggests us that not for application, still it might be of interests to look at $f(r)$-quasi spherical symmetry. That
is, bound the number of pieces one should cut an $r$-sphere in Cayley graph so that the harmonic measure ratio is say bounded by 2 in each piece.
Assumption \ref{Assumption1} might be true for some general $f(r)$-quasi spherically symmetric transitive infinite graphs.\\

From the proof of Theorem \ref{thm1}, what we really need is the following: {\it
There is a positive constant $c_6$ such that for any $n\in\mathbb{N},j\leq r_n$ and
any $A\subseteq \partial B_G(u,j)$ with $\mu_j(A)\geq \lambda_1(G_n)/2,$
$$\sum\limits_{y\in A}b_j(y,u,n)\geq c_6.\eqno{(4.2)}$$}
If (4.2) holds, then similarly to Theorem \ref{thm1}, we can verify that
without Assumption \ref{Assumption1}, Theorem \ref{thm1} does still hold.
To prove Conjecture \ref{conj1} in the uniformly nonamenable setting, the following problem should be studied.
\begin{problem}\label{pr3}
Let $G$ be a nonamenable transitive infinite connected graph. Is there a positive constant $c$ such that for any $n\in\mathbb{N}$ and
any $A\subseteq \partial B_G(u,n)$ with $\mu_n(A)\geq \lambda_1(G)/2,$
$$\sum\limits_{y\in A}b_n(y,u)\geq c?$$
\end{problem}

Furthermore, we can propose the following asymptotically absolutely continuous problem:
\begin{problem}\label{pr4}
Let $G$ be a nonamenable transitive infinite connected graph and each $\nu_n$ the uniform probability on $\partial B_G(u,n).$
{\bf (i)} Is $\mu_n$ (resp. $b_n(\cdot,u)$) asymptotically absolutely continuous with respect to $\nu_n$ in the sense that for any
$A_n\subseteq \partial B_G(u,n)$ with $\lim\limits_{n\rightarrow\infty}\nu_n(A_n)=0,$ we have
$$\lim\limits_{n\rightarrow\infty}\mu_n(A_n)=0\ \left(resp.\ \lim\limits_{n\rightarrow\infty}\sum\limits_{y\in A_n}b_n(y,u)=0\right)?$$
{\bf (ii)} Is $\mu_n$ asymptotically absolutely continuous with respect to $b_n(\cdot,u)$ in the sense that for any
$A_n\subseteq \partial B_G(u,n)$ with $\lim\limits_{n\rightarrow\infty}\sum\limits_{y\in A_n}b_n(y,u)=0,$ we have
$\lim\limits_{n\rightarrow\infty}\mu_n(A_n)=0?$
\end{problem}
Affirmative answer to Problem \ref{pr4}(i) or Problem \ref{pr4}(ii) can be used to prove Conjecture \ref{conj1} in the uniformly nonamenable
case similarly to Theorem \ref{thm1}.\\

Notice $\mathbb{Z}$ is an amenable 2-regular spherically symmetric infinite transitive graph and all examples for quasi-spherically symmetric
infinite transitive graphs in Section 2 are nonamenable; naturally the following problem arises:
\begin{problem}\label{pr5}
Are there amenable quasi-spherically symmetric infinite transitive connected graphs with degree at least 3?
\end{problem}

Finally, the following locality problem is very interesting in its own way. Let
$$p_u(G)= \inf\left\{p \in [0, 1] :\ \mathbb{P}_{G,p} (\exists ~\mbox{an unique infinite component})> 0\right\}.$$
Recall \cite{BS1996} conjectured that $p_c(G)<p_u(G)$ for any nonamenable quasi-transitive infinite
connected graph $G;$ and this conjecture holds in some cases (\cite{RLYP2013}).
\begin{problem}\label{pr6}
Let $\{G_n\}_{n=1}^{\infty}$ be a sequence of transitive infinite connected graphs converging locally to a transitive infinite
connected graph $G.$ {\bf (i)} Does $p_u(G_n) \rightarrow p_u(G)$ when each $G_n$ and $G$ are nonamenable?
{\bf (ii)} Under what conditions, $\theta_{G_n}(p)\rightarrow \theta_G(p)$ for any $p\in (0,1)$ or $\theta_{G_n}(\cdot)\rightarrow \theta_G(\cdot)$
in the Skorohod (resp. uniform) topology for C\`{a}dl\`{a}g (resp. continuous) functions on $[0,1]$?
\end{problem}

\noindent{\bf Acknowledgements.} Song He is indebted to Professors R. Lyons and V. Sidoravicius for useful discussions during ICM2014, Seoul.
The authors are very grateful to Professor G. Kozma for pointing out several mistakes in their manuscript, and to Professors I. Benjamini and R. Lyons
for helpful comments and suggestions. These improve greatly the quality of the paper.

\end{document}